\documentclass[12pt]{amsart}
\usepackage{amsmath}

\usepackage{amsmath,amssymb,amsthm,enumerate,color,graphicx}
\usepackage{tikz-cd}
\usepackage[colorlinks=true,citecolor=magenta,urlcolor=blue]{hyperref}

\usepackage[capitalise]{cleveref}
\crefname{theorem}{Theorem}{Theorems}
\crefname{lemma}{Lemma}{Lemmas}
\crefname{corollary}{Corollary}{Corollaries}
\crefname{proposition}{Proposition}{Propositions}
\crefname{conjecture}{Conjecture}{Conjectures}
\crefname{question}{Question}{Questions}
\crefname{definition}{Definition}{Definitions}
\crefname{example}{Example}{Examples}
\crefname{remark}{Remark}{Remarks}
\crefname{question}{Question}{Questions}
\crefname{enumi}{}{}
\crefname{equation}{}{}

\newtheorem{theorem}{Theorem}[section]
\newtheorem{corollary}[theorem]{Corollary}
\newtheorem{lemma}[theorem]{Lemma}
\newtheorem{proposition}[theorem]{Proposition}

\theoremstyle{definition}

\newtheorem{example}[theorem]{Example}

\theoremstyle{remark}
\newtheorem{remark}[theorem]{Remark}

\newcommand{\C}{\mathbb C}

\newcommand{\Z}{\mathbb Z}

\newcommand{\OO}{\mathcal O}

\newcommand{\ord}{\operatorname{ord}}

\newcommand{\beq}{\begin{equation}}
\newcommand{\eeq}{\end{equation}}

\makeatletter
\newsavebox\myboxA
\newsavebox\myboxB
\newlength\mylenA

\newcommand*\xoverline[2][0.9]{%
    \sbox{\myboxA}{$\m@th#2$}%
    \setbox\myboxB\null
    \ht\myboxB=\ht\myboxA%
    \dp\myboxB=\dp\myboxA%
    \wd\myboxB=#1\wd\myboxA
    \sbox\myboxB{$\m@th\overline{\copy\myboxB}$}
    \setlength\mylenA{\the\wd\myboxA}
    \addtolength\mylenA{-\the\wd\myboxB}%
    \ifdim\wd\myboxB<\wd\myboxA%
       \rlap{\hskip 0.5\mylenA\usebox\myboxB}{\usebox\myboxA}%
    \else
        \hskip -0.5\mylenA\rlap{\usebox\myboxA}{\hskip 0.5\mylenA\usebox\myboxB}%
    \fi}
		
\makeatother

\def\sideremark#1{\ifvmode\leavevmode\fi\vadjust{\vbox to0pt{\vss
 \hbox to 0pt{\hskip\hsize\hskip1em
 \vbox{\hsize2.7cm\tiny\raggedright\pretolerance10000
  \noindent #1\hfill}\hss}\vbox to8pt{\vfil}\vss}}}

\begin{document}
\title[]{Strongly Pseudoconvex Domains with Finite Dimensional $p$-Bergman Spaces}
\author{Shreedhar Bhat}
\address{Manipal Institute of Technology Bengaluru, Manipal Academy of Higher Education, Manipal, Karnataka, India}
\email{shreedhar.bhat@manipal.edu} 
\author{Soumya Ganguly}
\address{Department of Mathematics, Rutgers University, New Jersey}
\email{soumya.ganguly@rutgers.edu}
\author{Achinta Kumar Nandi}
\address{Department of Mathematics, University of California San Diego, La Jolla, CA}
\email{anandi@ucsd.edu} 

\date{\today}

\begin{abstract}
We show that for every $n\ge2$, $p>2$, and positive integer $m$, there exists a finite-volume pseudoconvex domain $D\subset\mathbb C^n$ with real-analytic strongly pseudoconvex boundary whose $p$-Bergman space $A^p(D)$ has dimension $m$. These domains admit smooth strictly plurisubharmonic exhaustions and satisfy $H^\infty(D)=\mathbb C$ and $\dim A^s(D)=\infty$ for every $0<s\le2$. Moreover, for every finite $s\ge p$, we have $\dim A^s(D)=\lceil pm/s\rceil$. As an application, we show that the Bergman projection fails to map $\mathcal{C}_c^\infty(D)$ into $A^p(D)$ and admits no bounded extension to $L^s(D)$ or $L^{s/(s-1)}(D)$ for any finite $s\ge p$.
\end{abstract}


\maketitle

\section{Introduction and Main Results}

Let $\Omega\subset\mathbb{C}^n$ be a domain. For
$0<p<\infty$, the $p$-Bergman space of $\Omega$ is
defined by
\[
A^p(\Omega)=\left\{
f\in\mathcal{O}(\Omega):
\int_\Omega |f|^p\,dV_{2n}<\infty
\right\},
\]
where $dV_{2n}$ denotes Euclidean Lebesgue measure.
For $p=2$, this is the classical Bergman space, whose
Hilbert space structure gives rise to the Bergman kernel
and, under suitable nondegeneracy conditions, the
Bergman metric; see~\cite{Bergman}.

A basic question is how the geometry of $\Omega$ and
the exponent $p$ determine the dimension of
$A^p(\Omega)$. Every bounded domain has an
infinite-dimensional $p$-Bergman space, since all
holomorphic polynomials are integrable there.
On the other hand, the submean inequality gives
$A^p(\mathbb{C}^n)=\{0\}$. Thus, domains with
positive finite-dimensional $p$-Bergman spaces
must be unbounded. On such domains, integrability
is influenced by both the boundary geometry
and the behavior of holomorphic functions at infinity.

The dimension problem has a long history.
In one complex variable, Carleson
\cite[Chapter VI]{Carleson}
characterized the nontriviality of the Bergman space
in terms of the logarithmic capacity of the complement
for domains with compact complement. The same
characterization holds for arbitrary planar domains;
see~\cite{GHH}. Wiegerinck~\cite{Wiegerinck}
proved in 1984 that
$\dim A^2(\Omega)\in\{0,\infty\}$
for every planar domain $\Omega$. In the same
paper, he constructed nonpseudoconvex Reinhardt
domains in $\mathbb{C}^2$ whose Bergman spaces
have any prescribed positive finite dimension.
This naturally leads to the question whether
the zero-or-infinity dichotomy holds for
pseudoconvex domains in higher dimensions,
commonly known as Wiegerinck's problem;
see~\cite{Jucha,PZBalanced}.

The dichotomy was established for several important classes of
pseudoconvex domains. Engli\v{s}~\cite{Englis} proved it for Reinhardt
domains, and Pflug--Zwonek~\cite{PZBalanced} for balanced domains in
$\mathbb{C}^2$. Jucha~\cite{Jucha} proved the dichotomy for complete
Hartogs domains over $\mathbb{C}$ defined by subharmonic weights,
characterizing their Bergman spaces through properties of the
defining potential. Boudreaux~\cite{Boudreaux} obtained further
sufficient conditions for infinite dimensionality for Hartogs
domains with higher-dimensional bases. These results make Hartogs
domains a natural setting for the dimension problem: expansion in
the fiber variables reduces global integrability to a family of
weighted problems on the base; see~\cite{Jucha,Boudreaux,BN}.

Finite-dimensional Bergman spaces also arise in complex geometry.
Huang--Li~\cite{HL} showed that positive constant holomorphic
sectional curvature of a well-defined Bergman metric forces the
Bergman space to be finite-dimensional. Hayashida--Kamimoto~\cite{HK} constructed such
domains in every dimension $n\ge2$ with constant
holomorphic sectional curvature $2$.
Bhat--Ganguly--Nandi--Xiao~\cite[Theorem 1.1]{BGNX}
subsequently constructed domains whose Bergman
metrics have constant holomorphic sectional curvature
$2/m$ for every positive integer $m$ and every $n\ge2$.
More recently, Yuan~\cite{Yuan} disproved the general pseudoconvex
dichotomy (i.e., Wiegerinck's problem) by constructing a pseudoconvex domain in $\mathbb{C}^2$
with a one-dimensional Bergman space, and pseudoconvex domains in
every $\mathbb{C}^n$, $n\geq3$, with arbitrary prescribed positive
finite Bergman dimension. These developments sharpen the question of which additional geometric hypotheses force infinite dimensionality, and how the answer changes with the integrability
exponent.

Recent developments in the theory of \(p\)-Bergman spaces have brought renewed attention to the dependence of their analytic and geometric properties on the integrability exponent \(p\). Chen-Zhang \cite{CZ} developed a general $p$-Bergman theory on bounded domains, including extremal kernels and reproducing formulas, while Chakrabarti-Edholm \cite{CE} constructed projections using monomial bases on Reinhardt domains. Bhat-Nandi \cite{BN} studied the dimension problem for $p$-Bergman spaces, obtaining geometric and potential-theoretic criteria through $L^p$ extension theorems in the range $1\leq p<2$. Carleson \cite{Carleson} characterized the nontriviality of Bergman spaces on planar domains in terms of the logarithmic capacity of their complements. The planar case also exhibits a distinction between different values of $p$. Positive finite-dimensional $p$-Bergman
spaces exist for $1\leq p<2$, whereas for $p\geq2$ the dimension is either zero or infinite (see \cite[Theorems 2.1 and 2.2]{BN}). Thus examples with positive finite $p$-Bergman dimension for $p>2$ necessarily require at least two complex variables.

\subsection{Main Results} In the present paper, we address this problem by constructing, for each prescribed positive integer $m$, a pseudoconvex domain whose $p$-Bergman space has dimension exactly $m$, for $p>2$. Our first main result is as follows.
\begin{theorem}\label{th:Hartogs}
For integers $n\ge2$ and $m\ge1$ and each $p>2$, there is a finite-volume pseudoconvex Hartogs domain $\Omega\subset\mathbb{C}^n$, with a smooth strictly plurisubharmonic exhaustion, such that its base coordinate $z$ satisfies
$$
A^p(\Omega)=\operatorname{span}_\mathbb{C}\{1,z,\ldots,z^{m-1}\}, \qquad H^\infty(\Omega)=\mathbb{C}.
$$
Furthermore, $\dim A^s(\Omega)=\infty$ for $0< s \leq 2$ and, for every finite $s \geq p$,
$$
A^s(\Omega)=\operatorname{span}_\mathbb{C} \{z^j:j\ge0,\ sj<pm\}, \qquad \dim A^s(\Omega)=\left\lceil\frac{pm}{s}\right\rceil.
$$ 
\end{theorem}
Motivated by Jucha's construction in \cite{Jucha}, we use logarithmic singularities, integrability conditions and infinite behavior of the subharmonic weight to gain control over the elements in the $p$-Bergman space. Suppose $D_\psi (X) \subset \mathbb{C}^n$ is a Hartogs domain with $n-1$ dimensional fibers with a planar domain $X$ as its base. From proposition $3.13$ in \cite{BN}, it follows that $F(z,w) =\displaystyle{ \sum_{\nu \in \mathbb{N}^{n-1}_0} f_\nu(z) w^\nu}$ with $f_\nu$ being in some weighted $p$-Bergman space, for all $F \in A^p(D_\psi(X))$. First by arranging logarithmic singularities of $\psi$ near the origin, we ensure that the local $L^p$ integrability force $f_\nu \equiv 0$ for all $|\nu| \geq 1$. In other words, using these logarithmic singularity arrangement we eliminate $w$ dependence of $F$. Next, using the infinite behavior of $\psi$ we ensure that any $F \in A^p(D_\psi(X))$ is a polynomial in $z$ of degree at most $m-1$. Finally adding a small $\log (1 + |z|^2)$ type term strict subharmonicity can be retained. See Proposition~\ref{prop:hartogs} in Section~\ref{sec2}.

The construction of the above Hartogs domains then can be pushed to construct unbounded domains with globally strictly pseudoconvex real analytic boundary. The Hartogs domains constructed above have exceptional boundary points where the boundary need not be smooth. Imitating the idea of a \textit{rational shear} due to Forstneri{\v c}-Wold in \cite{FW09} and \cite{FW13}, we use a holomorphic change of coordinates to eliminate the said weakly pseudoconvex points. The change of coordinates is constructed to have a constant Jacobian in order to preserve the unweighted $p$-Bergman spaces. Our second main result concerns construction of strictly pseudoconvex domains with real analytic boundary with prescribed finite dimensional $p$-Bergman spaces.

\begin{theorem}\label{th: strongly pseudoconvex}
Let $n \ge 2$ and $m \ge 1$ be integers, and let $p > 2$. There is a finite-volume pseudoconvex domain $D=D_{n,p,m} \subset \mathbb{C}^n$ whose boundary is everywhere real analytic and strongly pseudoconvex, such that
$$
\dim A^p(D)=m, \qquad H^\infty(D)=\mathbb{C}, \qquad \dim A^s(D)=\infty\quad(0<s\le2).
$$
The domain admits a smooth strictly plurisubharmonic exhaustion. Its first coordinate $u$ is nowhere zero, and for every finite
$s\ge p$,
\begin{equation}\label{eq:global-spaces}
A^s(D)=\operatorname{span}_{\mathbb{C}} \{u^{-j}:j \in \mathbb{Z}_{\ge0},\ sj<pm\}, \qquad \dim A^s(D)=\left\lceil \frac{pm}{s} \right\rceil.
\end{equation}  
\end{theorem}

As an application, we examine the $L^p$-regularity of the Bergman projection. We show that if $A^2(\Omega)$ is infinite-dimensional and $A^p(\Omega)$ is finite-dimensional, then the Bergman projection does not map $C_c^\infty(\Omega)$ into $A^p(\Omega)$. Moreover, the Bergman kernel $K_\Omega(\cdot,a)$ fails to belong to $A^p(\Omega)$ outside a proper analytic subset of $\Omega$. Applied to our examples, these results show that the Bergman projection is unbounded on $L^s$ and $L^{s/(s-1)}$ for every finite $s\geq p$, even though the domains have everywhere real-analytic strongly pseudoconvex boundary.

The restriction $p>2$ is essential to our construction. Following Jucha \cite{Jucha} and Bhat--Nandi \cite{BN}, we use logarithmic singularities to eliminate all terms involving the fiber variables. This requires an integer
$N>2/(p-2)$, and the argument does not work for $p=2$. In fact, our domains admit infinitely many linearly independent square-integrable holomorphic functions with simple poles at selected punctures. Hence $A^2$ is infinite-dimensional, as is $A^s$ for every $0<s\leq2$, since the domains have finite volume. This is also consistent with the result of Gallagher-Harz-Herbort
\cite{GHH} that a pseudoconvex domain with a smooth finite-type boundary point has an infinite-dimensional Bergman space.

\subsection*{Acknowledgement} The authors are grateful to Ming Xiao for valuable suggestions. The third author would like to thank Ji{\v r}{\' i} Lebl and Sean Curry for their inspiration and encouragement.

\subsection*{AI declaration}  The authors used OpenAI's ChatGPT to identify and study background literature, assist with grammar and sentence construction, and help assess the correctness and consistency of mathematical arguments. The authors take full responsibility for the content of this manuscript.

\section{The Hartogs Domains}\label{sec2}
Our ideas in the construction of the Hartogs domains with finite dimensional $p$-Bergman spaces are inspired from Jucha's construction in \cite{Jucha}. For complete Hartogs domains over planar bases, Jucha \cite{Jucha} used logarithmic singularities and integrability conditions to study the dimension of the Bergman space. In particular, his Proposition $4.5$ relates the existence of square-integrable holomorphic functions to the vanishing orders imposed by the logarithmic singularities of the defining weight. Related results for Hartogs domains and weighted spaces were obtained by Boudreaux \cite{Boudreaux} and Borichev-Le-Youssfi \cite{borichev2021dimension}. We use this relation in the opposite direction. Instead of starting with a Hartogs domain and determining its Bergman space, we construct a subharmonic weight for which the associated weighted spaces have prescribed dimensions.

We begin by briefly discussing the idea behind the construction of the Hartogs domain. Our goal is to construct a Hartogs domain $\Omega= \{(z,w): \|w\| < e^{-\psi_\varepsilon(z)}\}$ such that $A^p(\Omega) = \operatorname{span}_{\mathbb{C}} \{1, z, \dots, z^{m-1}\}$. To achieve this, using the subharmonic weight $\psi_\epsilon$ we must
\begin{itemize}
    \item eliminate every function that depends on $w$, and,
    \item ensure that the remaining functions are polynomials in $z$ with degree at most $m-1$.
\end{itemize}
Since any $F \in A^p(\Omega)$ is of the form $F(z,w) = \sum_\nu f_\nu (z) w^\nu$, to eliminate functions depending on $w$,  we need to ensure that $f_\nu \equiv 0$ whenever $|\nu| \geq 1$, . By making the $L^p$-integrability condition sufficiently restrictive near a point, one can force a holomorphic function to vanish to a prescribed order at that point (see Proposition \ref{prop:local-int}). However a single logarithmic singularity would not ensure $f_\nu \equiv 0$, we shall therefore need to use infinitely many such conditions. Our idea is to arrange infinitely many small collections of logarithmic singularities approaching the origin, each forcing more zeros than poles for every sufficiently large $|\nu|$. It follows from Jensen's formula that such non-zero coefficients must grow too rapidly in order to be $L^p$-integrable. Thus, we shall show that using appropriate local integrability conditions we can ensure that $F \in A^p(\Omega)$ is independent of $w$.

The remaining part is then to ensure that $F$ is a polynomial in $z$ of degree at most $m-1$. Clearly this requires controlling the growth of $F$ at infinity which can be done through the behavior of $\psi_\epsilon$ at infinity. Finally to force strict subharmonicity of $\psi_\epsilon$ we add a small $\log (1+|z|^2)$ term. Accordingly, we seek a family of subharmonic weights,
strictly subharmonic when $\varepsilon>0$,
of the following form:
$$
\psi_\varepsilon(z) = \sum_{\ell = 1}^\infty U_\ell (z) + \beta_\varepsilon \sum_j \log |z-b_j| + \varepsilon \log (1+|z|^2), 
$$
where the first term is used to remove $w$-dependence through local $L^p$ integrability, the middle term is for controlling the infinite behavior, and the last term is to ensure strict subharmonicity.

To this end, we begin the construction of the subharmonic weight by appropriately using logarithmic singularities. We denote the planar area measure by $dA$. For a nonzero meromorphic function of one variable, $\ord_a f$ is its integer order, positive at a zero and negative at a pole. Later, for a nonzero holomorphic germ $f$ at $a \in \mathbb{C}^n$, we use the same notation for its total vanishing order,
$$
\ord_a f=\min \bigl\{|\alpha|:\partial^\alpha f(a) \neq 0 \bigr\}.
$$
We begin with an elementary criterion for weighted integrability near an isolated puncture.
\begin{proposition}\label{prop:local-int}
Let $1 \le p < \infty$, $\lambda \in \mathbb{R}$, and $f \in \mathcal{O}(\{0<|z|<r_0\})$. If
\begin{equation}\label{eq:local-integral}
 \int_{0<|z|<r_0}|f(z)|^p|z|^{-p\lambda}\,dA(z) <\infty,
\end{equation}
then either $f=0$ or $f$ extends meromorphically across $0$ and
\begin{equation}\label{eq:local-order}
 \mathrm{ord}_0 f>\lambda-\frac2p.
\end{equation}
Conversely, a meromorphic function satisfying \eqref{eq:local-order} is integrable as in \eqref{eq:local-integral} on a sufficiently small punctured disk.
\end{proposition}

\begin{proof}
Writing $f(z) = \sum_{j \in \mathbb{Z}} a_j z^j$, using H{\" o}lder's inequality first note that
$$
a_j r^j = \frac{1}{2\pi} \int_0^{2\pi} f(re^{i\theta}) e^{-ij\theta} \,d\theta, \text{    and }\qquad |a_j|^pr^{pj} \leq \frac{1}{2\pi} \int_0^{2\pi} |f(re^{i \theta})|^p \,d \theta.
$$
Now whenever $a_j \neq 0$, it follows that
$$
\int_0^{r_0} |a_j|^p r^{pj - p\lambda +1} \,dr \leq \int_{\{0< |z| < r_0\}} |f(z)|^p |z|^{-p\lambda} \,dA(z).
$$
Therefore $pj - p\lambda + 1 > -1$ or $j > \lambda - \frac{2}{p}$. This proves the meromorphic extension and \eqref{eq:local-order}.

Conversely, if $\ord_0 f$ is $j$ then $f$ can be written as $z^j g(z)$ near $0$ for some holomorphic $g$ with $g(0) \neq 0$. Thus \eqref{eq:local-integral} holds if $j > \lambda - \frac{2}{p}$ in a smaller annulus where $g$ is non-vanishing.
\end{proof}

Let $d\ge1$ be an integer, let $X\subset\mathbb C$ be open, and let $\psi:X\to\mathbb R$ be continuous. Define
$$
D_{\psi} (X) := \{(z,w) \in X \times \mathbb{C}^d: \|w\| < e^{-\psi(z)} \}
$$
We denote $|\nu| = \sum_{j=1}^d \nu_j$ for a multi-index $\nu \in \mathbb{N}^d_{0}$, where $\mathbb{N}_0 := \mathbb{N} \cup \{0\}$. It follows from Proposition $3.13$ in \cite{BN} that every $F \in A^p(D_\psi(X))$ can be expressed in the form
$$
F(z,w) = \sum_{\nu \in \mathbb{N}^d_0} f_\nu(z) w^\nu,
$$
where $f_\nu \in \mathcal{O}(X)$ for all $\nu \in \mathbb{N}_0^d$ and $1 \leq p < \infty$. Moreover, as
$$
f_\nu(z) w^\nu = \frac{1}{(2\pi)^d} \int_{[0,2\pi]^d} F(z, e^{i\theta_1} w_1, \dots, e^{i\theta_d} w_d) e^{-i \nu. \theta} \,d\theta
$$
it follows that
\begin{equation}\label{eq:coefficient-integral}
\int_X |f_\nu(z)|^p e^{-(p|\nu| + 2d) \psi(z)} \,dA(z) < \infty,
\qquad \text{  for each } \nu.
\end{equation}

The next lemma is a form of Jensen's formula on the annuli. The construction below will force the coefficients of positive fiber degree to have more zeros than poles in each of infinitely many small annuli accumulating at the origin. We use Jensen's formula to show that a nonzero meromorphic function satisfying these conditions must have rapidly growing circular means near the origin. The following lemma allows the zeros and poles within each annulus to occur at different radii, provided the annuli are sufficiently narrow.

\begin{lemma}\label{lem:jensen}
Let $f$ be a nonzero meromorphic function on $0<|z|<R$, and let $\mathcal{A}_\ell =\{e^{-b_\ell}\le |z|\le e^{-a_\ell}\}$ be pairwise disjoint annuli tending to the origin. Suppose that every pole of $f$ sufficiently close to the origin belongs to one of these annuli. Let $N_\ell$ and $P_\ell$ denote, respectively, the numbers of zeros and poles of $f$ in $\mathcal A_\ell$, counted with multiplicity. Assume that $N_\ell-P_\ell \geq 1$ for all sufficiently large $\ell$, and that
$$
\sum_{\ell}(b_\ell-a_\ell)P_\ell<\infty.
$$
Then, for every $L>0$, there is a constant $C_L$ such that
$$
\frac{1}{2\pi} \int_0^{2\pi} \log|f(re^{i\theta})|\,d\theta \geq L\log \frac{1}{r} -C_L
$$
for all sufficiently small $r>0$.
\end{lemma}

\begin{proof}
For $0 < r= e^{-t} < R$, set
\[
M_f(t)=\frac{1}{2\pi}\int_0^{2\pi} \log|f(e^{-t+i\theta})|\,d\theta.
\]
Choose $t_0$ sufficiently large so that every pole of $f$ in $0<|z|<e^{-t_0}$ belongs to one of the annuli $\mathcal A_\ell$, and $N_\ell-P_\ell\ge1$ for every annulus lying in this punctured disc. We may also assume that $f$ has no zeros or poles on $|z|=e^{-t_0}$.

By Jensen's formula, there are constants $A,B\in\mathbb R$ such that, for $t \ge t_0$,
\begin{equation}\label{eq:jensen-representation}
M_f(t)=A+ Bt + \sum_{\{\zeta:\,t_0<-\log|\zeta|\le t\}} \ord_{\zeta}(f) \bigl(t+\log|\zeta|\bigr).
\end{equation}
Here zeros are counted with positive multiplicities and poles with negative multiplicities. The sum is finite for each fixed $t$, since $f$ is meromorphic on the punctured disc.

We estimate the contribution from each annulus $\mathcal A_\ell$. Recall that $\mathcal{A}_\ell =\{e^{-b_\ell}\le |z|\le e^{-a_\ell}\}.$ Suppose first that $t \ge b_\ell$, so that the circle $|z|=e^{-t}$ lies on or inside the inner boundary circle of $\mathcal A_\ell$. Every zero $\zeta$ in this annulus satisfies
$$
t + \log|\zeta| \ge t-b_\ell,
$$
whereas every pole satisfies
$$
t+\log|\zeta|\le t-a_\ell.
$$
Consequently, the total contribution from $\mathcal A_\ell$ to \eqref{eq:jensen-representation} is at least
\begin{align*}
N_\ell(t-b_\ell)-P_\ell(t-a_\ell) &=(N_\ell-P_\ell)(t-b_\ell) -(b_\ell-a_\ell)P_\ell\\
&\ge (t-b_\ell) -(b_\ell-a_\ell)P_\ell.
\end{align*}

Next, suppose that $a_\ell\le t<b_\ell$. In this case, only some of the zeros and poles in $\mathcal A_\ell$ occur in the sum \eqref{eq:jensen-representation}. The zeros contribute nonnegative terms, while the total
negative contribution from the poles is bounded below by
$$
-(b_\ell-a_\ell)P_\ell.
$$
Finally, any zeros outside the annuli contribute nonnegative terms. By our choice of $t_0$, there are no poles outside the annuli in $0<|z|<e^{-t_0}$. Combining these estimates, we obtain
$$
M_f(t)\ge A+Bt+ \sum_{\{\ell:\,a_\ell>t_0\\b_\ell\le t\}} (t-b_\ell) -\sum_\ell(b_\ell-a_\ell)P_\ell.
$$
By hypothesis, $C:=\sum_\ell(b_\ell-a_\ell)P_\ell<\infty.$ Hence
\begin{equation}\label{eq:jensen-lower-bound}
M_f(t)\ge A+Bt+ \sum_{\{\ell:\,a_\ell>t_0,\\
b_\ell\le t\}}(t-b_\ell)-C.
\end{equation}
Now let $L>0$. Choose an integer $J$ large enough that $B+J>L.$ Since there are infinitely many annuli approaching the origin, we may choose $J$ of them lying in $0 < |z| < e^{-t_0}$. For all sufficiently large $t$, these $J$ annuli occur in the sum in \eqref{eq:jensen-lower-bound}.
Every other term in that sum is nonnegative. Therefore, if $b_{\ell_1},\ldots,b_{\ell_J}$ are the right endpoints of the selected intervals, then
$$
M_f(t) \ge A+Bt+\sum_{j=1}^J(t-b_{\ell_j})-C =(B+J)t+A-C-\sum_{j=1}^Jb_{\ell_j} \ge Lt-C_L
$$
for a suitable constant $C_L$.

This proves the lemma.
\end{proof}

We now begin to construct the logarithmic weights. Fix $n \geq 2, p > 2, m \geq 1$ throughout this section, set $d=n-1$ and
$$
N = 1 + \left\lfloor \frac{\frac{2}{p}}{1-\frac{2}{p}} \right\rfloor = 1+  \left\lfloor \frac{2}{p-2} \right\rfloor.
$$
In particular $0 < \frac{2}{p} < 1$ and $\frac{2}{pN} < 1- \frac{2}{p}$. 

If $\psi (z) = \alpha \log |z-a| + O(1)$ near an isolated deleted point $a$, Proposition \ref{prop:local-int} shows that
$$
\int |f|^p e^{-(2(n-1)+pk) \psi} \,dA < \infty \implies \ord_a f > \alpha \left( \frac{2}{p}(n-1) + k \right) - \frac{2}{p},
$$
and since $\ord_a f$ is an integer, it follows that $\ord_a f \geq 1 + \left \lfloor \alpha \left( \frac{2}{p}(n-1) + k \right) - \frac{2}{p} \right \rfloor$. 

Define $R_k (\alpha):= 1 + \left\lfloor \alpha \left( \frac{2}{p}(n-1) + k \right) - \frac{2}{p} \right\rfloor$, for all integers $k \geq 0$. We shall arrange the coefficients of the logarithmic terms so that these order bounds force every coefficient of positive fiber degree to vanish.

For each $j \geq 1$ and $k \geq 1$, define
$$
a_j = \frac{2}{2(n-1) + pj}, \quad D_{kj} = R_k(a_j) + N R_k \left(-\frac{a_j}{N}\right).
$$
For each $j \geq 1$, we use an $(N+1)$-tuple consisting of one $a_j$ and $N$-many $-\frac{a_j}{N}$ coefficients. Clearly their sum is zero and $(n-1)a_j < 1$. We record some properties of $D_{kj}$ in the following lemma.
\begin{lemma}\label{lem:integers}
The integers $D_{kj}$ satisfy
$$
D_{kj} = 0 \quad \text{ for } j > k, \quad D_{kk} = 1, \quad R_0(a_j) = R_0(-\frac{a_j}{N}) = 0.
$$
Further, there are $M_j \in \mathbb{N}$ such that $\sum_{j = 1}^\ell M_j D_{k j} \geq 1$, for all $\ell \geq k \geq 1$.
\end{lemma}

\begin{proof}
To show that $D_{kj} = 0$ for all $j > k$, first note that
$$
N R_k(-\frac{a_j}{N}) = N + N \left\lfloor -\frac{2}{pN} \frac{\left(\frac{2}{p}(n-1) + k\right)}{\left(\frac{2}{p}(n-1) + j\right)} - \frac{2}{p} \right\rfloor, \text{  and   } \frac{\left(\frac{2}{p}(n-1) + k\right)}{\left(\frac{2}{p}(n-1) + j\right)} < 1.
$$
Using $\frac{2}{pN} < 1- \frac{2}{p}$ from above, for all $j > k$, it follows that
$$
0 > -\frac{2}{pN} \frac{\left(\frac{2}{p}(n-1) + k\right)}{\left(\frac{2}{p}(n-1) + j\right)} - \frac{2}{p} > - \frac{2}{p} - \frac{2}{pN} > -1.
$$
Thus $N R_k( -\frac{a_j}{N}) = N - N = 0$. Similarly,
$$
R_k(a_j) = 1 + \left\lfloor \frac{2}{p} \frac{\left(\frac{2}{p}(n-1) + k \right)}{\left(\frac{2}{p}(n-1) + j \right)} - \frac{2}{p} \right\rfloor = 1-1 =0.
$$
Thus, $D_{kj} = 0$ for all $j > k$.

To prove $D_{kk} = 1$, it is straightforward to see that $R_k(a_k) = 1$ and $N R_k \left( -\frac{a_k}{N} \right) = N-N = 0$. Also, $R_0(a_j) = 0$ and $R_0(-\frac{a_j}{N}) = 0$ is also immediate from the definition of $R_k$ above. 

As $D_{kk} = 1$ and $D_{kj} = 0$ for all $j > k$, it follows that $\sum_1^\ell M_j D_{kj} = M_k + \sum_1^{k-1} M_j D_{kj}$. Set, $M_1 = 1$ and $M_k = 1 + \sum_1^{k-1} M_j \max \{-D_{kj}, 0\}$. It follows that
$$
\sum_1^\ell M_j D_{kj} = 1 + \sum_1^{k-1} M_j (D_{kj} + \max \{ 0, -D_{kj} \}) \geq 1.
$$

This completes the proof.
\end{proof}

For each $\ell \geq 1$, we consider $M_j$ copies of the $N+1$-tuple $\left(a_j, -\frac{a_j}{N}, \dots, -\frac{a_j}{N} \right)$, for all $1 \leq j \leq \ell$. The total number of coefficients $J_\ell$ and the sum of their absolute values $V_\ell$ are
$$
J_\ell = (N+1) \sum_{j=1}^\ell M_j, \quad V_\ell = 2\sum_{j=1}^\ell a_j M_j.
$$
Set $T_\ell = J_\ell + V_\ell$.

We now choose where to place the logarithmic singularities. Since the coefficients in each collection sum to zero, placing the corresponding points sufficiently close together makes their logarithmic potential small away from those points. The following lemma also gives the integral estimate needed to control the volume of the Hartogs domain.

Fix $\ell\ge1$ and choose distinct real numbers
$$
\xi_{i,s,j}\in(-1/2,1/2), \qquad 1\le i\le\ell,\quad 1\le s\le M_i,\quad 0\le j\le N.
$$
For $r\in\mathbb C$ and $\delta>0$, define
$$
U_\delta(z)= \sum_{i=1}^{\ell}\sum_{s=1}^{M_i}a_i \left(\log|z-r-\delta\xi_{i,s,0}|
-\frac1N\sum_{j=1}^{N} \log|z-r-\delta\xi_{i,s,j}| \right).
$$
\begin{lemma}\label{lem:log-sum}
For every fixed $\rho > 0$,
\begin{equation}\label{eq:log-sum-limits}
\sup_{|z-r|\ge\rho}|U_\delta(z)| \to 0, \qquad \int_{|z-r| < \rho} e^{-2(n-1)U_\delta(z)} \,dA(z) \to \pi\rho^2
\end{equation}
as $\delta \downarrow 0$.
\end{lemma}

\begin{proof}
For fixed $\rho > 0$ and for all $|z-r| \geq \rho$, we have
$$
U_\delta (z) = \sum_{i=1}^\ell \sum_{s = 1}^{M_i} a_i \left( \log \left|1- \frac{\delta \xi_{i,s,0}}{z-r} \right| - \frac{1}{N} \sum_{j=1}^N \log \left| 1 - \frac{\delta \xi_{i,s,j}}{z-r} \right| \right).
$$
For $|z-r|\ge\rho$ and $0<\delta<\rho$,
we have
\[
\left|\frac{\delta\xi_{i,s,j}}{z-r}\right|
\le\frac{\delta}{2\rho}<\frac12.
\]
Since $\log|1-\eta|=O(|\eta|)$ uniformly for
$|\eta|\le1/2$, and the sum defining
$U_\delta$ has finitely many terms, we obtain
\[
\sup_{|z-r|\ge\rho}|U_\delta(z)|
=O\left(\frac{\delta}{\rho}\right)
\longrightarrow0.
\]

Say, $U(\zeta) = \sum_{i=1}^\ell \sum_{s=1}^{M_i} a_i \left( \log |\zeta - \xi_{i,s,0}| - \frac{1}{N} \sum_{j=1}^N \log |\zeta - \xi_{i,s,j}| \right)$. It follows that $U(\zeta) = O(|\zeta|^{-1})$ as $|\zeta| \to \infty$. Note that near $\xi_{i,s,0}$ the function $e^{-2(n-1)U(\zeta)}$ is comparable to $|\zeta - \xi_{i,s,0}|^{-2(n-1) a_i}$, and therefore $e^{-2(n-1)U(\zeta)} = 1 + O(|\zeta|^{-1})$ for large enough $|\zeta|$. Writing
$$
\int_{|\zeta| < R} e^{-2(n-1) U(\zeta)} \,dA(\zeta) = \int_{|\zeta| < R_0} e^{-2(n-1) U(\zeta)} \,dA(\zeta) + \int_{R_0 \leq |\zeta| < R} e^{-2(n-1) U(\zeta)} \,dA(\zeta),
$$
it follows that $\int_{|\zeta| < R_0} e^{-2(n-1) U(\zeta)} \,dA(\zeta) = O(1)$, and $\int_{R_0 \leq |\zeta| < R} e^{-2(n-1) U(\zeta)} \,dA(\zeta) = \pi R^2 + O(R)$. Here $R_0$ is chosen so that all logarithmic singularities are contained in $|\zeta| < R_0$.

For the second assertion, set $\zeta=(z-r)/\delta$.
Since $U_\delta(r+\delta\zeta)=U(\zeta)$, we obtain
\[
\begin{aligned}
\int_{|z-r|<\rho}e^{-2(n-1)U_\delta(z)}\,dA(z)
&=\delta^2
  \int_{|\zeta|<\rho/\delta}
  e^{-2(n-1)U(\zeta)}\,dA(\zeta)\\
&=\delta^2\left(
  \pi\frac{\rho^2}{\delta^2}
  +O\left(\frac{\rho}{\delta}\right)+O(1)
  \right)\\
&=\pi\rho^2+O(\rho\delta+\delta^2).
\end{aligned}
\]
The asserted limit follows as $\delta\downarrow0$.
\end{proof}

We now apply the lemma to construct infinitely many collections of singularities approaching the origin. Using these logarithmic singularities we ensure that any $L^p$-integrable holomorphic function on the Hartogs domain must be independent of $w$. 

Recall that $T_\ell = J_\ell + V_\ell$. Let
\begin{equation}\label{eq:disks}
r_\ell=2^{-\ell}, \qquad \rho_\ell= \frac{r_\ell2^{-\ell-10}}{1+T_\ell},
\qquad B_\ell=\{z\in\mathbb C:|z-r_\ell|<\rho_\ell\}.
\end{equation}
The disks $B_\ell$ are pairwise disjoint and lie in $\{|z| <1 \}$. Moreover,
\begin{equation}\label{eq:width-sum}
\sum_{\ell=1}^{\infty} T_\ell\frac{\rho_\ell}{r_\ell}<\infty,
\end{equation}
since
$$
T_\ell\frac{\rho_\ell}{r_\ell} =\frac{T_\ell}{1+T_\ell}2^{-\ell-10} \le 2^{-\ell-10}.
$$
For each $\ell$, choose distinct real numbers
$$
\xi_{\ell,i,s,j}\in(-1/2,1/2), \qquad 1\le i\le\ell,\quad 1\le s\le M_i,\quad 0\le j\le N.
$$
By Lemma \ref{lem:log-sum}, we may choose $0<\delta_\ell<\rho_\ell/4$ such that, with
\begin{equation}\label{eq:locations}
A_{\ell,i,s,j} = r_\ell+\delta_\ell\xi_{\ell,i,s,j}
\end{equation}
and
$$
U_\ell(z) =\sum_{i=1}^{\ell}\sum_{s=1}^{M_i}a_i \bigg(\log|z-A_{\ell,i,s,0}| -\frac{1}{N} \sum_{j=1}^{N} \log|z-A_{\ell,i,s,j}| \bigg),
$$
we have
\begin{equation}\label{eq:local-estimates}
\sup_{z\notin B_\ell}|U_\ell(z)| \le 2^{-\ell}, \quad \int_{B_\ell}e^{-2dU_\ell}\,dA \le 2\pi\rho_\ell^2.
\end{equation}

All the points $A_{\ell,i,s,j}$ are distinct and real, and their only accumulation point is $0$. By Lemma \ref{lem:log-sum}, we may choose $0 < \delta_\ell < \frac{\rho_\ell}{4}$ such that \eqref{eq:local-estimates} holds.

\subsection{The Defining Functions} The singularities constructed above will be used to eliminate all positive fiber degrees. We now introduce additional logarithmic terms to control the degree of the remaining holomorphic functions. For this, we want the defining function to satisfy $e^{-2d\psi_\varepsilon(z)} \asymp |z|^{-2-pm}$ as $|z| \to \infty.$ Indeed, this weight makes $z^k$ integrable near infinity precisely when $k<m$. Equivalently, we require
$$
\psi_\varepsilon(z)=\tau\log|z|+O(1), \qquad \tau=\frac{1+pm/2}{d}.
$$
We obtain this growth by adding logarithmic singularities at finitely many points $b_1,\ldots,b_{\gamma_\varepsilon}$, where $\tau_\varepsilon=\tau-2\varepsilon$ and $\gamma_\varepsilon=\lceil2(n-1)\tau_\varepsilon\rceil$. The total coefficient of these logarithms must be $\tau_\varepsilon$, since the term $\varepsilon\log(1+|z|^2)$ contributes $2\varepsilon\log|z|$ at infinity.

To ensure local integrability, we divide this coefficient equally among sufficiently many points. More precisely, for $0 \le \varepsilon < \frac{\tau}{2},$ set
$$
\beta_\varepsilon =\frac{\tau_\varepsilon}{\gamma_\varepsilon},
$$
and choose $b_j=j+2, 1\le j\le \gamma_\varepsilon.$ Then $\gamma_\varepsilon\beta_\varepsilon =\tau-2\varepsilon,$ and $ 0 < (n-1) \beta_\varepsilon \le \frac{1}{2}.$ In particular, the weight $|z-b_j|^{-2(n-1)\beta_\varepsilon}$ is locally integrable at each $b_j$.

We can now define the punctured base and the subharmonic function $\psi_\varepsilon$.

Define
$$
K_\varepsilon =\{0,b_1,\ldots,b_{\gamma_\varepsilon}\} \cup \bigl\{ A_{\ell,i,s,j}: \ell\ge1,\ 1\le i\le\ell,\ 1\le s\le M_i,\ 0\le j\le N \bigr\},
$$
and $X_\varepsilon=\mathbb C\setminus K_\varepsilon.$ Finally, let
\begin{equation}\label{eq:background}
h_\varepsilon(z) =\beta_\varepsilon \sum_{j=1}^{\gamma_\varepsilon}\log|z-b_j|
+\varepsilon\log(1+|z|^2),
\end{equation}
and define
\begin{equation}\label{eq:psi}
\psi_\varepsilon(z) = h_\varepsilon(z)+\sum_{\ell=1}^{\infty}U_\ell(z), \qquad z \in X_\varepsilon.
\end{equation}
The set $K_\varepsilon$ is compact and countable. Since it lies on the real axis, its complement $X_\varepsilon$ is connected.

By \eqref{eq:local-estimates}, the series $\sum_{\ell\ge1}U_\ell$ converges uniformly on compact subsets of $X_\varepsilon$. Its sum is harmonic there. Thus $\psi_\varepsilon$ is real analytic and
\begin{equation}\label{eq:laplacian}
(\psi_\varepsilon)_{z\bar z} = \frac{\varepsilon}{(1+|z|^2)^2}.
\end{equation}
In particular, $\psi_\varepsilon$ is harmonic when $\varepsilon=0$ and strictly subharmonic when $\varepsilon>0$.

At $A_{\ell,i,s,0}$, we have
$$
\psi_\varepsilon(z) =a_i\log|z-A_{\ell,i,s,0}|+O(1),
$$
whereas, for $1\le j\le N$,
$$
\psi_\varepsilon(z) =-\frac{a_i}{N} \log|z-A_{\ell,i,s,j}|+O(1).
$$
Similarly, at each $b_j$,
$$
\psi_\varepsilon(z) =\beta_\varepsilon\log|z-b_j|+O(1).
$$

The estimates \eqref{eq:local-estimates} imply
\begin{align}
|\psi_\varepsilon-h_\varepsilon|
&\le1
&&\text{outside }\bigcup_\ell B_\ell,
\label{eq:outside}\\
|\psi_\varepsilon-h_\varepsilon-U_\ell|
&\le1
&&\text{on }B_\ell\setminus K_\varepsilon.
\label{eq:inside}
\end{align}

Since all the disks $B_\ell$ lie in the unit disk,
the sum $\sum_\ell U_\ell$ is bounded at infinity.
Also,
$$
h_\varepsilon(z) = \bigl(\gamma_\varepsilon\beta_\varepsilon
+2\varepsilon\bigr)\log|z|+O(1) =\tau\log|z|+O(1).
$$
Consequently,
\begin{equation}\label{eq:infinity}
\psi_\varepsilon(z)=\tau\log|z|+O(1), \qquad 2(n-1)\tau=2+pm.
\end{equation}

Finally, consider the annuli
\begin{equation}\label{eq:annuli}
G_\ell = \left\{\frac35r_\ell<|z|<\frac45r_\ell \right\}.
\end{equation}
These annuli are disjoint from every $B_j$. Since $h_\varepsilon$ is bounded near the origin, \eqref{eq:outside} gives a constant $C_\varepsilon$ such that
\begin{equation}\label{eq:annular-bound}
|\psi_\varepsilon|\le C_\varepsilon \qquad \text{on  }\bigcup_{\ell\ge1}G_\ell.
\end{equation}

\subsection{The Hartogs Domain} Using the function $\psi_\varepsilon$ defined in \eqref{eq:psi}, we set
$$
\Omega_\varepsilon = D_{\psi_\varepsilon}(X_\varepsilon) := \{(z,w) \in X_\varepsilon \times \mathbb{C}^{n-1} : \|w\| < e^{-\psi_\varepsilon(z)} \}.
$$
The following proposition summarizes the above construction.
\begin{proposition}\label{prop:hartogs} For every $0 \le \varepsilon < \tau/2$, the domain $\Omega_{\varepsilon}$ is connected, pseudoconvex, and of finite volume, and admits a smooth strictly plurisubharmonic exhaustion. For every finite $s\ge p$,
\begin{equation}\label{eq:hartogs-spaces}
 A^s(\Omega_{\varepsilon}) =\operatorname{span}_{\mathbb{C}}\{z^j:j \in \mathbb{Z}_{\ge0},\ sj<pm\}, \quad \dim A^s(\Omega_{\varepsilon}) =\left\lceil\frac{pm}{s}\right\rceil.
\end{equation}
Moreover, $H^\infty(\Omega_{\varepsilon})=\mathbb{C}$ and $\dim A^s(\Omega_{\varepsilon}) = \infty$, for $0 <s \le 2$. If $\varepsilon > 0$, the boundary over $X_\varepsilon$ is real analytic and strongly pseudoconvex.
\end{proposition}

Note that, when $\varepsilon > 0$ the portion of $\partial \Omega_{\varepsilon}$ contained in $X_\varepsilon \times \mathbb{C}^{n-1}$ is real analytic and strongly pseudoconvex. However, the full boundary need not be smooth, because it also contains points lying over
$K_\varepsilon$. We shall split the proof of Proposition \ref{prop:hartogs} above in several smaller observations.

We shall start with the geometric assertions.
\begin{proposition}\label{prop:geometry}
The domain $\Omega_{\varepsilon}$ is connected, has finite volume, and admits a smooth strictly plurisubharmonic exhaustion. In particular, it is pseudoconvex.
\end{proposition}

\begin{proof}
To prove connectedness, note that $X_\varepsilon$ is connected and contracting each fiber of $\Omega_\varepsilon$ we get a retract to $X_\varepsilon$.

To prove finiteness of the volume of $\Omega_\varepsilon$, first note that
$$
\text{vol}(\Omega_\varepsilon) = \text{vol}(\mathbb{B}^{n-1}) \int_{X_\varepsilon} e^{-2(n-1) \psi_\varepsilon} \,dA.
$$
Since $(n-1) \beta_\varepsilon \leq \frac{1}{2}$ near each $b_j$ the term $e^{-2(n-1) h_\varepsilon}$ is integrable, and it is integrable near infinity because $h_\varepsilon = \tau \log |z| + O(1)$ with $(n-1)\tau > 1$. By \eqref{eq:outside}, on
$X_\varepsilon\setminus\bigcup_{\ell\ge1}B_\ell$
we have $e^{-2d\psi_\varepsilon}
\le e^{2d}e^{-2dh_\varepsilon}.$
The function $e^{-2dh_\varepsilon}$ is locally
integrable near each $b_j$, since
$2d\beta_\varepsilon\le1<2$.
Moreover, at infinity it is comparable to
$|z|^{-2d\tau}$, where
$2d\tau=2+pm>2$.
It follows that the volume integral over
$X_\varepsilon\setminus\bigcup_{\ell\ge1}B_\ell$
is finite. On the union of $B_\ell$, the function $h_\varepsilon$ is uniformly
bounded. Combining equations \eqref{eq:inside} and \eqref{eq:local-estimates}, we obtain 
$$
\int_{B_\ell\setminus K_\varepsilon} e^{-2d\psi_\varepsilon}\,dA
\le C\int_{B_\ell}e^{-2dU_\ell}\,dA \le 2\pi C\rho_\ell^2.
$$
The right-hand side is summable. This proves that $\text{vol}(\Omega_\varepsilon)$ is finite.

To prove pseudoconvexity, we start by enumerating $K_\varepsilon \setminus \{0\} = \{\zeta_j: j \geq 1\}$ and define
$$
\chi (z) = |z|^2+|z|^{-2} + \sum_{j \geq 1}2^{-j}|z-\zeta_j|^{-2}.
$$
On compact subsets of $X_\varepsilon$, the function $\chi$ and all its derivatives converge uniformly. Also, each pole term is subharmonic away from its pole, and the first term is strictly subharmonic. The function tends to infinity as $z$ approaches $K_\varepsilon$ or infinity. Thus $\chi$ is a smooth strictly subharmonic exhaustion of $X_\varepsilon$.

We now extend this exhaustion to the fibers as follows. Set, $S(z,w) = \|w\|^2 e^{2\psi_\varepsilon(z)}$. The function $\log \|w\|^2+2\psi_\varepsilon(z)$ is
plurisubharmonic where $w\ne0$. More directly, the Levi
form of $S=\|w\|^2e^{2\psi_\varepsilon}$ satisfies
\[
\mathcal L_S(v,\xi)
=e^{2\psi_\varepsilon}
\left(
\|\xi+2(\partial\psi_\varepsilon)(v)w\|^2
+2\|w\|^2\mathcal L_{\psi_\varepsilon}(v)
\right)\ge0.
\]
This formula also holds at $w=0$, where
$\mathcal L_S(0,\xi)=e^{2\psi_\varepsilon}\|\xi\|^2$.
Thus $S$ is plurisubharmonic and its Levi form is
positive definite on nonzero vertical vectors. Therefore
\begin{equation}\label{eq:exhaustion}
\chi(z)+\frac{S(z,w)}{1-S(z,w)}
\end{equation}
is strictly plurisubharmonic as $t/(1-t)$ is increasing and convex on
$[0,1)$, while $\chi$ is strict in every nonzero base direction. Its sublevel sets keep $z$ in a compact subset of $X_\varepsilon$ and $S$ bounded away from one, so it is an exhaustion.
\end{proof}

For $k\geq0$, and $f \in \mathcal{O}(X_\varepsilon)$, we set
$q_k=k+\frac{2(n-1)}{p}$, and define
$$
I_{k, \varepsilon}(f) := \int_{X_\varepsilon} |f(z)|^p e^{-(pk + 2(n-1))\psi_{\varepsilon}(z)} \,dA(z)
$$
\begin{proposition}\label{prop:positive-degree}
For each integer $k \geq 1$, if $f \in \mathcal{O}(X_\varepsilon)$ satisfies $I_{k,\varepsilon}(f) < \infty$ then $f \equiv 0$.
\end{proposition}

\begin{proof}
The idea is that the singularities accumulating at the origin force $f$ to have more zeros than poles in each of infinitely many small disks. If $f$ is non-zero, by Jensen's formula it follows that the circular means grow too rapidly near the origin for the weighted integral to be finite.

Fix $k \ge 1$ and suppose that $f\not\equiv0$ satisfies $I_{k,\varepsilon}(f)<\infty$. By Proposition \ref{prop:local-int}, $f$ extends meromorphically across each of the points $A_{\ell,i,s,j}$, with
\[
\ord_{A_{\ell,i,s,0}}f\ge R_k(a_i), \quad \text{ and }
\ord_{A_{\ell,i,s,j}}f \ge R_k(-a_i/N),\qquad 1\le j\le N.
\]
Since $D_{ki} = R_k (a_i) + N R_k \left(-\frac{a_i}{N} \right)$, for every $\ell \ge k$, it follows that
$$
\sum_{i=1}^{\ell}\sum_{s=1}^{M_i} \sum_{j=0}^{N}\ord_{A_{\ell,i,s,j}}f \ge \sum_{i=1}^{\ell}M_i D_{ki} \ge 1,
$$
where the last inequality follows from Lemma \ref{lem:integers}. Thus, in every disk $B_\ell$ with $\ell\ge k$, the number of zeros of $f$ exceeds the number of poles, counted with multiplicities, by at least one.

We next verify the hypotheses of Lemma \ref{lem:jensen}. Since $R_k(a_i)\ge0$, poles can occur only at the points $A_{\ell,i,s,j}$ with $1\le j\le N$. Moreover,
$$
\max\{0,-R_k(-a_i/N)\} \le \frac{q_ka_i}{N}+1.
$$
Hence the total pole multiplicity $P_{k,\ell}$ in $B_\ell$ satisfies
\begin{equation}\label{eq:pole-bound}
P_{k,\ell}  \le \sum_{i=1}^{\ell} M_iN\left(\frac{q_ka_i}{N}+1\right) \le q_kV_\ell+J_\ell \le(q_k+1)T_\ell.
\end{equation}

Recall that $B_\ell=\{|z-r_\ell|<\rho_\ell\}$. To apply Lemma~\ref{lem:jensen}, consider the pairwise disjoint annuli $\{r_\ell-\rho_\ell\le |z|\le r_\ell+\rho_\ell\}$. Each contains $B_\ell$, and since all poles sufficiently close to the origin lie in these disks, passing to the annuli introduces no additional poles. Thus the excess of zeros over poles remains at least one for $\ell\ge k$. As $z$ ranges over $B_\ell$, the quantity $t=-\log|z|$ lies between $-\log(r_\ell+\rho_\ell)$ and $-\log(r_\ell-\rho_\ell).$ The difference between these endpoints is
$$
\log\frac{r_\ell+\rho_\ell}{r_\ell-\rho_\ell} \le\frac{4\rho_\ell}{r_\ell}.
$$
Therefore, by \eqref{eq:width-sum} and \eqref{eq:pole-bound}, we obtain
\[
\sum_{\ell=1}^{\infty} P_{k,\ell} \log\frac{r_\ell+\rho_\ell}{r_\ell-\rho_\ell}
<\infty.
\]
All poles sufficiently close to the origin lie in these disks. We may therefore apply
Lemma \ref{lem:jensen}, which gives, for every $L>0$,
\begin{equation}\label{eq:positive-mean-growth}
M_f(t):= \frac1{2\pi}\int_0^{2\pi} \log|f(e^{-t+i\theta})|\,d\theta \ge Lt-C_L
\end{equation}
for all sufficiently large $t$.

We now show that this contradicts the integrability of $f$. Recall that the annuli $G_\ell= \left\{\frac35r_\ell<|z|<\frac45r_\ell\right\}$ are disjoint from all the disks $B_j$. By \eqref{eq:annular-bound}, the weight $e^{-pq_k\psi_\varepsilon}$ is bounded below by a positive constant on their union. By Jensen's inequality, we have
\[
\frac1{2\pi}\int_0^{2\pi} |f(e^{-t+i\theta})|^p\,d\theta \ge e^{pM_f(t)}.
\]
Using polar coordinates and then substituting $r=e^{-t}$, we obtain
\[
I_{k,\varepsilon}(f) \ge C_k\sum_{\ell=1}^{\infty} \int_{-\log(4r_\ell/5)}^{-\log(3r_\ell/5)}
e^{pM_f(t)-2t}\,dt
\]
for some constant $C_k>0$. Choose $L>2/p$. By \eqref{eq:positive-mean-growth}, for all sufficiently large $\ell$,
\[
\begin{aligned}
I_{k,\varepsilon}(f) &\ge C_{k,L}\sum_{\ell\ge\ell_0} \int_{-\log(4r_\ell/5)}^{-\log(3r_\ell/5)} e^{(pL-2)t}\,dt.
\end{aligned}
\]
Each interval of integration has length $\log(4/3)$, and its endpoints tend to infinity as $\ell\to\infty$. Since $pL-2>0$, the sum diverges.

This contradicts $I_{k,\varepsilon}(f)<\infty$. Therefore $f\equiv0$.    
\end{proof}

Thus as $F \in A^p(\Omega_\varepsilon)$ are of the form $F(z,w) = \sum_{\nu \in \mathbb{N}^{n-1}_0} f_\nu(z) w^\nu$ with $I_{|\nu|, \varepsilon} (f) < \infty$, by Proposition \ref{prop:positive-degree} every function in $A^p(\Omega_\varepsilon)$ is independent of $w$. We now treat the remaining degree zero terms using infinite behavior of the subharmonic weight $\psi_\varepsilon$.

\begin{proposition}\label{prop:degree-zero}
The functions $f\in\OO(X_\varepsilon)$ with $I_{0,\varepsilon}(f)<\infty$ are exactly the polynomials of degree at most $m-1$.
\end{proposition}

\begin{proof} At every point $A_{\ell,i,s,j}$, the lower bound for the order is zero, since $R_0(a_i)=R_0(-a_i/N)=0$. At a point $b_j$, it is also zero, since
$$
-1<- \frac{2}{p} \leq q_0\beta_\varepsilon- \frac{2}{p} =\frac{2}{p}((n-1) \beta_\varepsilon - 1) < 0.
$$
Therefore $f$ extends holomorphically through every isolated point of $K_\varepsilon$, and hence to $\C\setminus\{0\}$.

The Laurent expansion near $0$ is $f(z)=\sum_{j\in\Z}a_jz^j$. The following inequality is well known
\begin{equation}\label{eq:circle}
 |a_j|^pr^{pj}\le \frac1{2\pi}\int_0^{2\pi}|f(re^{i\theta})|^p\,d\theta.
\end{equation}
For every $s \in \mathbb{N}$, using \eqref{eq:circle} and \eqref{eq:annular-bound}, we obtain 
\begin{align*}
 I_{0,\varepsilon}(f) \ge C|a_{-s}|^p\sum_{\ell\gg1}\int_{3r_\ell/5}^{4r_\ell/5}r^{1-ps} \,dr =C_{p,s}|a_{-s}|^p\sum_{\ell\gg1}r_\ell^{2-ps}.
\end{align*}
Since $p > 2$, the last sum diverges. All negative Laurent coefficients vanish, so $f$ is entire.

At infinity, \eqref{eq:infinity} shows that $e^{-2d\psi_\varepsilon(z)}$ is comparable to $|z|^{-2d\tau}$. Applying \eqref{eq:circle} to the Taylor series of $f$, any nonzero coefficient $a_j$ must satisfy
$\int_R^\infty r^{pj-2d\tau+1}\,dr<\infty,$ equivalently $pj<2d\tau-2=pm.$ Thus $j<m$.

Conversely, any polynomial of degree at most $m-1$ is bounded on bounded subsets of the base. Proposition \ref{prop:geometry} forces integrability over those subsets, and the same strict inequality proves integrability at infinity.
\end{proof}

\begin{proposition}\label{prop:spaces}
For all $s \geq p > 2$, $A^s(\Omega_{\varepsilon}) = \operatorname{span}_{\mathbb{C}} \{z^j: j \geq 0, sj < pm\}$. Moreover, $H^\infty(\Omega_\varepsilon) = \mathbb{C}$ and
$\dim A^s(\Omega_{\varepsilon})=\infty$ for every $0<s\le2$. The domain $\Omega_\varepsilon$ is not biholomorphic to a balanced domain.
\end{proposition}

\begin{proof}
From Proposition \ref{prop:degree-zero} it follows that $A^p(\Omega_\varepsilon) = \operatorname{span}_{\mathbb{C}} \{1, z, \dots, z^{m-1}\}$.

Since $\Omega_\varepsilon$ has finite volume (see Proposition \ref{prop:geometry}), it follows that $A^s (\Omega_\varepsilon) \subset A^p(\Omega_\varepsilon)$ for all $s \geq p$. A polynomial of degree $j$ is in $A^s$ precisely when its integral at infinity converges, namely when $sj<2(n-1)\tau-2 = pm$. Its integral over a bounded part of the base is finite by Proposition~\ref{prop:geometry}. This proves the formula for $A^s$. The number of integers $j\ge0$ with $j< \frac{pm}{s}$ is $\lceil \frac{pm}{s} \rceil$.

Next, note that every bounded holomorphic function belongs to $A^p$, so it must be a polynomial in $z$. Since $X_\varepsilon$ contains the exterior of a disk, only a constant polynomial can be bounded. Hence $H^\infty=\mathbb{C}$.

At each point $A_{\ell,i,s,j}$ with $1\le j\le N$, we have $\psi_\varepsilon(z) =-\frac{a_i}{N}\log|z-A_{\ell,i,s,j}|+O(1).$ So, the function $F_{\ell,i,s,j}(z,w) = (z-A_{\ell,i,s,j})^{-1}$ belongs to $A^2 (\Omega_{\varepsilon})$. Indeed, its squared norm near $A_{\ell,i,s,j}$ is comparable to
\[
\int_{|z-A_{\ell,i,s,j}|<r} |z-A_{\ell,i,s,j}|^{-2+2(n-1)a_i/N}\,dA(z) <\infty.
\]
Away from $A_{\ell,i,s,j}$, the function $F_{\ell,i,s,j}$ is bounded, so its square
integrability follows from the finite volume of $\Omega_{\varepsilon}$. Since there are infinitely many distinct points $A_{\ell,i,s,j}$ with $1\le j\le N$, the
corresponding functions $F_{\ell,i,s,j}$ are linearly independent. Hence $\dim A^2(\Omega_{\varepsilon})=\infty.$ Since $\Omega_{\varepsilon}$ has finite volume it follows that $A^2(\Omega_{\varepsilon}) \subset A^s (\Omega_{\varepsilon}),$ for $0 < s < 2$.

Finally, $\Omega_\varepsilon$ is not even homeomorphic to any balanced domain because the non-contractible $X_\varepsilon \subset \mathbb{C}$ is a deformation retract of $\Omega_\varepsilon$, whereas every balanced domain is star-shaped about the origin and therefore is contractible.

This completes the proof.
\end{proof}

\begin{proposition}\label{prop:levi}
If $\varepsilon>0$, every point of
\[
\Gamma_\varepsilon =\{(z,w):z\in X_\varepsilon,\ \|w\|=e^{-\psi_\varepsilon(z)}\}
\]
is a real-analytic strongly pseudoconvex boundary point of $\Omega_{\varepsilon}$. The set $\Gamma_\varepsilon$ is nonempty and relatively open in the boundary.
\end{proposition}

\begin{proof}
At a point of $\Gamma_\varepsilon$, $w \ne 0$, and $\rho(z,w) = \log\|w\| + \psi_\varepsilon(z)$ is a real-analytic local defining function with nonzero differential. For a vector $(v,\xi)\in\C\times\C^d$, its Levi form is
$$
\mathcal L_\rho(v,\xi) =\frac{\varepsilon}{(1+|z|^2)^2}|v|^2 +\frac1{2\|w\|^2} \left(\|\xi\|^2 - \frac{|\langle\xi,w\rangle|^2}{\|w\|^2}\right).
$$
The first term is positive if $v\ne0$. If $v=0$ and $(v,\xi)$ is complex tangent, then $\langle\xi,w\rangle=0$, so the second term is positive unless $\xi=0$. This proves strong pseudoconvexity.

Finally, $\Gamma_\varepsilon = \partial\Omega_\varepsilon \cap (X_\varepsilon \times\mathbb{C}^{n-1}).$ Since $X_\varepsilon$ is open, $\Gamma_\varepsilon$ is relatively open in $\partial\Omega_\varepsilon$.
\end{proof}

\subsection{Proof of Proposition \ref{prop:hartogs}}
\begin{proof}
The assertions follow from Propositions~\ref{prop:geometry}, \ref{prop:spaces}, and \ref{prop:levi}.
\end{proof}

\subsection{Proof of Theorem \ref{th:Hartogs}}
\begin{proof} The proof follows from setting $\varepsilon=0$ in Proposition \ref{prop:hartogs}.
\end{proof}

\subsection{An Example} 
\begin{example}\label{ex:p3} We give an example of a pseudoconvex domain $\Omega \subset \mathbb{C}^2$ of finite volume such that
$$
A^3(\Omega)=\mathbb C, \qquad \dim A^2(\Omega)=\infty.
$$

We take $n=2$, $p=3$, and $m=1$ in the preceding construction, with $\varepsilon=0$. In this case, $N=3, a_i=\frac{2}{3i+2}.$ The first four rows and columns of the matrix $(D_{ki})$ are
\[
\begin{pmatrix}
1&0&0&0\\
-2&1&0&0\\
-2&1&1&0\\
-1&-2&1&1
\end{pmatrix},
\]
and the first four multiplicities are
\[
M_1=1,\qquad M_2=3,\qquad M_3=3,\qquad M_4=8.
\]
The remaining multiplicities are determined recursively by Lemma~\ref{lem:integers}.

For each $\ell\ge1$, choose distinct real points
\[
A_{\ell,i,s,j}, \qquad
1\le i\le\ell,\quad
1\le s\le M_i,\quad
0\le j\le3,
\]
sufficiently close to $2^{-\ell}$, as in the preceding construction. These points accumulate only at the origin.

Set $X=\mathbb C\setminus \left( \{0,3,4,5,6,7\} \cup \{A_{\ell,i,s,j}\} \right).$
Define $\psi:X \to \mathbb R$ by
\[
\begin{aligned}
\psi(z)
={}&\frac12\sum_{b=3}^{7}\log|z-b|
+\sum_{\ell=1}^{\infty}
\sum_{i=1}^{\ell}
\sum_{s=1}^{M_i}
\frac{2}{3i+2}
\left(
\log|z-A_{\ell,i,s,0}|
-\frac13\sum_{j=1}^{3}
\log|z-A_{\ell,i,s,j}|
\right).
\end{aligned}
\]
The points $A_{\ell,i,s,j}$ have been chosen so that the series converges locally uniformly on $X$. Thus $\psi$ is harmonic, and hence subharmonic, on $X$.

We now define
\[
\Omega= \left\{ (z,w)\in X\times\mathbb C: |w|<e^{-\psi(z)} \right\}.
\]
This is a pseudoconvex domain of finite volume.

The logarithmic singularities at the points $A_{\ell,i,s,j}$ eliminate every coefficient of positive degree in $w$, by Proposition \ref{prop:positive-degree}. So elements in $A^3(\Omega)$ depends only on $z$. The five additional logarithmic terms give $\psi(z)=\frac52\log|z|+O(1)$ as $|z|\to\infty.$ Hence $e^{-2\psi(z)}\asymp |z|^{-5}$ as $|z|\to\infty.$ In particular, the monomial $z^j$ is integrable to the third power near infinity precisely when $\int_1^\infty r^{3j-4}\,dr<\infty,$
or equivalently, $j<1$.

The preceding results therefore imply $A^3(\Omega)=\mathbb{C}$, and $A^s(\Omega)=\mathbb C$ for $3\le s<\infty,$ and $H^\infty(\Omega)=\mathbb C.$ On the other hand, $\dim A^2(\Omega)=\infty,$ and $\dim A^s(\Omega)=\infty$ for $0 < s \le 2.$
\end{example}

\section{Strongly Pseudoconvex Domains}
The Hartogs domain constructed above has the desired $p$-Bergman spaces, but its boundary is not strongly pseudoconvex at points lying over $K_\varepsilon = \mathbb{C} \setminus X_\varepsilon$. Following the idea of the rational shears of Wold and Forstneri{\v c}-Wold (see \cite{FW09,FW13}), we use a holomorphic change of coordinates to send these boundary points to infinity. We choose the map to have constant Jacobian so that it preserves the unweighted Bergman spaces. The resulting domain has an everywhere real-analytic strongly pseudoconvex boundary.

\begin{lemma}\label{lem:coordinates} Let $K \subset \mathbb{C}$ be compact and countable, with $0 \in K$, and suppose that every point of $K \setminus \{0\}$ is isolated. Set $X=\mathbb C\setminus K$. Let $R:X\to(0,\infty)$ be continuous and suppose that
\begin{equation}\label{eq:R-hypotheses}
\begin{aligned}
|z-a|R(z)&\longrightarrow0
&&\text{as }z\to a,\quad a\in K\setminus\{0\},\\
R(z)|z|^{-1} & \longrightarrow 0
&&\text{as }|z|\to\infty.
\end{aligned}
\end{equation}
Then there exists a holomorphic function $G$ on $X$
such that
\begin{equation}\label{eq:Phi}
\Phi(z,w_1,w') = \left(\frac1z,z^2(w_1+G(z)),w'\right)
\end{equation}
maps the domain $\mathcal B_R= \{(z,w)\in X\times\mathbb C^{n-1}: \|w\|<R(z)\}$
biholomorphically onto a domain $D \subset \mathbb{C}^n$, with complex Jacobian of $\Phi$ being $\det \Phi'=-1.$ Moreover, every boundary point of $D$ has a neighborhood on which $\Phi^{-1}$ is holomorphic and maps the boundary of $D$ onto the boundary of $\mathcal B_R$ over $X$. Here $w'\in\mathbb C^{n-2}$ is omitted when $n=2$.
\end{lemma}

\begin{proof}
We construct $G$ with simple poles at the points
of $K\setminus\{0\}$. Write
$K\setminus\{0\}=\{a_j\}$ and set
\[
G(z)=z+\sum_j\frac{2^{-j}}{z-a_j},
\]
where the sum is finite if $K\setminus\{0\}$
is finite. The series converges absolutely and uniformly on compact subsets of $X$. Since each $a_j$ is isolated,
\begin{equation}\label{eq:G-pole}
G(z)=\frac{2^{-j}}{z-a_j}+O(1)
\qquad(z\to a_j).
\end{equation}
Also, since $K$ is bounded,
\begin{equation}\label{eq:G-infinity}
G(z)=z+O(|z|^{-1})
\qquad(|z|\to\infty).
\end{equation}

The first coordinate of $\Phi$ determines $z$, and its inverse is $\Phi^{-1}(u,v,w') = \left(\frac1u,u^2v-G(1/u),w'\right).$ Thus $\Phi$ is biholomorphic onto its image. Moreover, $\det \Phi' =-\frac1{z^2}\cdot z^2=-1.$

We now show that no boundary point lying over $K$, or arising from $z \to \infty$, can have a finite image under $\Phi$. Suppose that $(z_k,w_k)\in\mathcal B_R$ and $\Phi(z_k,w_k)$ converges in $\mathbb C^n$. If $z_k\to0$, then $1/z_k$ is unbounded, which is impossible.

If $z_k\to a_j \ne 0$, then $|w_{1,k}+G(z_k)| \ge |G(z_k)|-R(z_k) \to  \infty,$
because $G$ has a simple pole at $a_j$ and $|z-a_j|R(z)\to0$. Since $a_j\ne0$, it follows that $|z_k^2(w_{1,k}+G(z_k))| \to \infty,$ again a contradiction. Finally, if $|z_k| \to \infty$, then $|w_{1,k}+G(z_k)| \ge |G(z_k)|-R(z_k) \to \infty,$ by \eqref{eq:R-hypotheses} and \eqref{eq:G-infinity}. Thus, $|z_k^2(w_{1,k}+G(z_k))| \to \infty,$ which is also impossible. It follows that every finite boundary point of $D$ has first coordinate $u \ne 0$, with $1/u\in X$. The inverse $\Phi^{-1}$ is holomorphic in a neighborhood of such a point. In this neighborhood, it identifies $D$ with $\mathcal{B}_R$. Since $R$ is continuous, it maps the boundary of $D$ onto
$$
\{(z,w):z\in X,\ \|w\|=R(z)\}.
$$
This proves the lemma.
\end{proof}

\subsection{Proof of Theorem \ref{th: strongly pseudoconvex}}

\begin{proof} Fix $0 < \varepsilon < \frac{\tau}{2}$ and consider the Hartogs
domain
$$
\Omega_\varepsilon =\left\{(z,w)\in X_\varepsilon\times\mathbb C^{n-1}:
\|w\|<e^{-\psi_\varepsilon(z)}\right\}.
$$
We first verify that $R = e^{-\psi_\varepsilon}$ satisfies the hypotheses of Lemma \ref{lem:coordinates}.

At each point $A_{\ell,i,s,0}$, we have $R(z)=O\left(|z-A_{\ell,i,s,0}|^{-a_i}\right),$ where $a_i=\frac{2}{pi+2(n-1)}<1.$ At the points $A_{\ell,i,s,j}$ with $1\le j\le N$, the function $R$ tends to zero. Similarly, at each $b_j$,
$$
R(z)=O(|z-b_j|^{-\beta_\varepsilon}), \qquad 0<\beta_\varepsilon\le\frac{1}{2(n-1)}<1.
$$
Therefore, $|z-a|R(z) \to 0$ as $z\to a$ at every point $a \in K_\varepsilon \setminus \{0\}$. Moreover, by \eqref{eq:infinity}, $R(z)\asymp |z|^{-\tau}$ near $|z|\to\infty,$ so $R(z)/|z|\to0$.

Thus we obtain a biholomorphism from Lemma \ref{lem:coordinates}
$$
\Phi(z,w_1,w') =\left(\frac1z,z^2(w_1+G(z)),w'\right)
$$
from $\Omega_\varepsilon$ onto a domain $D \subset \mathbb{C}^n$, with Jacobian determinant equal to $-1$.
Consequently, $\operatorname{vol}(D)
= \operatorname{vol}(\Omega_\varepsilon) <\infty.$
Therefore any holomorphic $f$ belongs to $A^s(D)$ iff $f \circ \Phi \in A^s (\Omega_\varepsilon)$. Writing $(u,v,w')=\Phi(z,w_1,w')$, we have $z=1/u$. Thus for every finite $s\ge p$, from Proposition \ref{prop:hartogs} it follows that $A^s(D) =\operatorname{span}_{\mathbb C} \{u^{-j}:j\ge0,\ sj<pm\},$ and hence $\dim A^s(D)=\left\lceil\frac{pm} {s}\right\rceil.$ The same proposition forces $H^\infty(D)=\mathbb{C},$
$\dim A^s(D)=\infty$ for all $0<s\le2$.

Since $\Omega_\varepsilon$ admits a smooth strictly plurisubharmonic exhaustion, composing it with $\Phi^{-1}$ gives such an exhaustion of $D$. In particular, $D$ is pseudoconvex.

Finally, Lemma \ref{lem:coordinates} shows that every boundary point of $D$ corresponds, under a local biholomorphism, to a point of $\Gamma_\varepsilon$.
By Proposition \ref{prop:levi}, $\Gamma_\varepsilon$ is real analytic and strongly pseudoconvex. These properties are preserved by biholomorphisms, so the entire boundary of $D$ is real analytic and strongly pseudoconvex.
\end{proof}

\begin{remark}\label{rem:intermediate-exponents}
The construction also yields information about $A^s(D)$ for $2<s<p$. By \eqref{eq:infinity} and the finite volume of $\Omega_\varepsilon$, we have $z^j\in A^s(\Omega_\varepsilon)$ whenever $sj<pm$. Since $\Phi$ preserves Lebesgue measure, it follows that
\[
\operatorname{span}_{\mathbb C}\{u^{-j}:j\ge0,\ sj<pm\}\subset A^s(D),\qquad
\dim A^s(D)\ge\left\lceil\frac{pm}{s}\right\rceil\ge m+1.
\]
In fact, $A^s(D)$ remains infinite-dimensional for certain $s>2$. Recall that $a_1=2/(p+2(n-1))$ and $N=1+\lfloor2/(p-2)\rfloor$, and set
\[
s_0=2+\frac{2(n-1)a_1}{N}
=2+\frac{4(n-1)}{N(p+2(n-1))}<p.
\]
At each point $A_{\ell,1,1,1}$, we have $\psi_\varepsilon(z)=-(a_1/N)\log|z-A_{\ell,1,1,1}|+O(1)$. Consequently, $(z-A_{\ell,1,1,1})^{-1}\in A^s(\Omega_\varepsilon)$ whenever $s<s_0$. Indeed, its integrability near the pole is equivalent to
\[
\int_0^\rho t^{1-s+2(n-1)a_1/N}\,dt<\infty,
\]
where $\rho>0$ is sufficiently small, while integrability elsewhere follows from the finite volume of $\Omega_\varepsilon$. Since these functions have distinct poles as $\ell$ varies, they are linearly independent. The volume-preserving biholomorphism $\Phi$ therefore gives $\dim A^s(D)=\infty$ for every $0<s<s_0$.
\end{remark}

\section{An application to the Bergman projection}\label{sec:projection-app}

He, Tang, and Tu \cite{HeTangTu} showed that the Bergman projection need not be bounded on $L^p$ for unbounded strongly pseudoconvex domains with real-analytic boundary. Related results for smooth unbounded worm domains were obtained by
Krantz, Monguzzi, Peloso, and Stoppato \cite{KMPS}. For our examples when $s \geq p$, the failure of $L^s$-boundedness follows from the finite-dimensionality of $A^p$ and the infinite-dimensionality of $A^2$. We first establish a general result.

Let $P_\Omega$ denote the orthogonal Bergman projection on $L^2(\Omega)$, and let
$K_\Omega(z,a)$ be the Bergman kernel.

\begin{proposition}\label{prop:projection}
Let $\Omega\subset\mathbb{C}^n$ be a domain such that $\dim A^2(\Omega)=\infty, \dim A^p(\Omega)<\infty$, for some $p > 0$. Then there exists a proper analytic
subset $E\subset\Omega$ such that, for every $a\in\Omega\setminus E$,
\[
K_\Omega(\cdot,a)\notin A^p(\Omega).
\]
Moreover, for every such $a$, there exists $\eta_a\in C_c^\infty(\Omega)$ satisfying
$P_\Omega\eta_a=K_\Omega(\cdot,a).$ In particular, $P_\Omega(C_c^\infty(\Omega))
\not\subset A^p(\Omega).$ If $p \geq 1$, the Bergman projection has no bounded
extension to $L^p(\Omega)$.
\end{proposition}

\begin{proof}
Since $A^2(\Omega)\cap A^p(\Omega)$ is finite-dimensional, we may choose a nonzero
function $g\in A^2(\Omega)$ orthogonal to this intersection. If $K_\Omega(\cdot,a) \in A^p(\Omega)$, then $g(a) =\langle g,K_\Omega(\cdot,a)\rangle=0.$ Thus $K_\Omega(\cdot,a)\notin A^p(\Omega)$ whenever $g(a) \ne 0$. So we can set
$E=\{a\in\Omega:g(a) = 0\}.$

Pick $a\in\Omega\setminus E$. Choose a nonnegative radially symmetric function
$\eta_a\in C_c^\infty(\Omega)$, supported in a ball centered at $a$, such that
\[
\int_\Omega\eta_a\,dV=1.
\]
For every $h\in A^2(\Omega)$, due to the mean-value
property $\langle h,\eta_a\rangle =h(a) =\langle h,K_\Omega(\cdot,a)\rangle.$
Hence $P_\Omega\eta_a=K_\Omega(\cdot,a) \notin A^p(\Omega).$ Finally, $\eta_a\in L^2(\Omega)\cap L^p(\Omega)$. If $p\geq1$, a bounded extension of $P_\Omega$
to $L^p(\Omega)$ would imply $P_\Omega\eta_a\in L^p(\Omega)$, a contradiction.
\end{proof}

\begin{corollary}\label{cor:projection-examples}
Let $D=D_{n,p,m}$ be a domain from Theorem \ref{th: strongly pseudoconvex}, where $p>2$. There exists $\eta\in C_c^\infty(D)$ such that $$P_D\eta\notin L^s(D)$$ whenever $p\leq s\leq\infty.$ Moreover, there is a proper analytic subset $E\subset D$ such that $K_D(\cdot,a) \notin L^s(D)$ for every $a\in D\setminus E$ and $p \leq s \leq \infty$. Consequently, for every finite $s\geq p$, the Bergman projection has no bounded extension to either $L^s(D)$ or $L^{s/(s-1)}(D)$.
\end{corollary}

\begin{proof}
By Theorem~\ref{th: strongly pseudoconvex},
$$
\dim A^2(D)=\infty, \qquad \dim A^p(D)=m, \qquad \mathrm{vol}(D)<\infty.
$$
Proposition~\ref{prop:projection} gives $\eta\in C_c^\infty(D)$ such that $P_D\eta \notin L^p(D)$. Since $D$ has finite volume, $L^s(D)\subset L^p(D)$ for all $p\leq s\leq\infty$. Thus $P_D\eta\notin L^s(D)$ for every $s\geq p$. The same argument applies to the Bergman kernel outside the proper analytic subset given by the proposition.

In particular, $P_D$ cannot be bounded on $L^s(D)$ for any finite $s\geq p$.
If it were bounded on $L^{\frac{s}{(s-1)}}(D)$, its self-adjointness on $L^2(D)$ and
$L^s$-duality would force $P_D\eta \in L^s(D)$, a contradiction.
\end{proof}

\begin{remark}
The corollary concerns global $L^s$-regularity
and does not imply a failure of Condition R.
In particular, failure of $L^s$-integrability
on an unbounded domain does not imply failure
of smooth extension at its finite boundary points.
\end{remark}

\bibliographystyle{alpha}
\bibliography{references}

\end{document}